\documentclass[preprint,12pt]{elsarticle}
\usepackage{amssymb}
\usepackage{lineno}
\usepackage[tbtags]{amsmath}
\usepackage{graphicx,floatrow,threeparttable}

\begin{document}
\begin{frontmatter}
\title{Hereditary effects of exponentially damped oscillators with past histories}
\author[rvt]{Jian Yuan \tnoteref{mytitlenote}\corref{cor}}
\address{College of Mathematic and Information Science, Shandong Institute of Business and Technology, Yantai 264005, P.R.China}
\cortext[cor]{Corresponding author. Email: yuanjianscar@gmail.com}
\author[rvt]{Guozhong Xiu\tnoteref{mytitlenote}}
\author{Bao Shi, Liying Wang}
\address{Institute of System Science and Mathematics, Naval Aeronautical University, Yantai 264001, P.R.China}
\tnotetext[mytitlenote]{These authors contributed equally to this work.}
\begin{abstract}
Hereditary effects of exponentially damped oscillators with past histories are considered in this paper. Nonviscously damped oscillators involve hereditary damping forces which depend on time-histories of vibrating motions via convolution integrals over exponentially decaying functions. As a result, this kind of oscillators are said to have memory. In this work, initialization for nonviscously damped oscillators is firstly proposed. Unlike the classical viscously damped ones, information of the past history of response velocity is necessary to fully determine the dynamic behaviors of nonviscously damped oscillators. Then, initialization response of exponentially damped oscillators is obtained to characterize the hereditary effects on the dynamic response. At last, stability of initialization response is proved and the hereditary effects are shown to gradually recede with increasing of time.
\end{abstract}

\begin{keyword}
Nonviscously damping model
\sep exponentially damped oscillators
\sep initialization problems
\sep hereditary effects
\end{keyword}

\end{frontmatter}

\section{Introduction}
Viscoelastic materials have seen broad applications in vibration control engineering due to their high damping capacity [1]. The modeling of constitutive relations of viscoelastic materials is a fundamental task for analysis and design of viscoelastically damped structures. However, it is also a challenging work because the mechanical behavior of viscoelastic materials is highly dependent on various factors such as time, temperature, the vibrating frequency and so on [2].

The integral constitutive models are derived based on the materials properties of stress relaxation and creep. The stress relaxation functions and creep functions are memory and hereditary kernels in the integral constitutive equations. They can be expressed by a series of exponential functions [3, 4], power-law functions [5], Mittage-Leffler functions [5], or other types of functions. The integral constitutive models are superior to the differential ones in many aspects: (i) the fading memory property can be characterized and time-history of loading acting on materials can be recorded; (ii) the stress relaxation functions or creep functions are easily and directly obtained via experiment data fitting; (iii) other factors such as temperature and the ageing affects can be conveniently considered and included in the model. The integral constitutive relations can be readily used as damping models in viscoelastically damped structures. Equations of motion of such systems are a set of coupled second-order integro-differential equations. The presence of the "integral" term makes the vibration analysis and control design more complicated than the classical ones. The integral type damping models may be also called nonviscously damping models and the corresponding oscillators are called nonviscously damped oscillators.

Researches on nonviscously damped oscillators are mainly concentrated on two types: one is the exponentially damped oscillators where the damping force are expressed by exponentially fading memory kernel; the other is the fractional-order oscillators where the viscoelastic relaxation functions are characterized by power-law functions or Mittage-Leffler functions. Adhikari and his colleagues have systematically investigated the structural dynamics with exponentially damped models, including dynamics of exponentially damped single-degree-of-freedom and multi-degree-of-freedom systems, identification and quantification of damping in [6-8]. They have discovered that, unlike the classical viscously damped oscillators, an exponentially damped single-degree-of-freedom oscillator has three eigenvalues. The complex conjugate pair of roots corresponds to the vibration motion and the third one corresponds to a purely dissipative motion which is always non-oscillatory in nature. The dynamics of exponentially damped oscillators is governed by both of the viscous damping factor and the nonviscous damping factor. For the dynamic analysis of multi-degree-of-freedom systems, they have developed a state-space approach using additional dissipation coordinates and the configuration space method. It was shown that the characteristic equation for an $N$ degree-of-freedom system is more than $N$ and the modes are divided into elastic modes and nonviscous modes. In [9], an analytical solution using modal superposition and two numerical solutions on the basis of finite element formulations are developed for the analysis of an exponentially damped solid rod. In [10], a method has been proposed to calculate eigensolution derivatives for the nonviscously damped systems. In [11], a state-space method has been proposed to identify modal and physical parameters for the nonviscously damped systems using recorded time-histories of the dynamic response. A closed-form approximation expression of the eigenvalues for non-viscous, non-proportional vibrating system has been derived in [12]. The forced vibration response to an arbitrary forced excitation has been obtained in [4].

Fractional-order oscillators are another type of nonviscously damped oscillators. They are under extensive investigations for the last three decades. The constitutive models involving fractional-order derivatives have been viewed more accurate and concise than other ones, because the relaxation functions of viscoelastic materials can be perfectly fitted from experimental data by using power-law functions or Mittage-Leffler functions [13, 14]. As a result, researches on fractional oscillators have been expected to be a promising work for structural dynamics analysis and control design. However, the fractional differential equations of motions are difficult to deal with due to the presence of the weakly singular kernels. Studies on dynamic responses of fractional-order oscillators have been reviewed in [15]. Asymptotically steady state behavior of fractional oscillators has been presented in [16, 17]. The criteria for the existence and the behavior of solutions on the basis of functional analytic approach have been obtained in [18-20], and the impulsive response function for the linear single-degree-of-freedom fractional oscillator has been derived. The asymptotically steady state response of fractional oscillators with more than one fractional derivatives have been analyzed in [21]. Based on the energy storage and dissipation properties of the Caputo fractional-order derivatives, the expression of mechanical energy in single-degree-of-freedom fractional oscillators has been determined and energy regeneration and dissipation during the vibratory motion have been obtained in [22]. Vibration controls have been designed using sliding mode control technique and adaptive control technique for single-degree-of-freedom fractional oscillators, multi-degree-of-freedom fractional oscillators, and fractional Duffing oscillators in [23, 24].

The damping forces in nonviscously damped oscillators depend on the past history of motion via convolution integrals over relaxation functions. As a result, the dynamics of the nonviscously damped oscillators is said to have memory. However, the memory and hereditary properties of fractional-order systems have been neglected and ambiguous for a long time [25]. Considering the memory effect and prehistory of fractional oscillators, the history effects and initialization problems for fractionally damped vibration equations have been proposed in [26-30]. The stability of initialization response has been proved based on the unit impulse response function and the Lyapunov stability theorem in [31].

This paper investigates the hereditary effects on the dynamics of the exponentially damped oscillators with past histories. We first declare that knowledge of the equations of motion, along with the initial displacement and velocity is insufficient to determine the dynamics behaviors. The initial conditions for the systems should also contain the past history of response velocity. Then we obtain the initialization response of exponentially damped oscillators, which characterizes the hereditary effects of the history of vibration on the dynamic response. At last, we prove that hereditary effects on the initialization response recede to zero with increasing of time.

\section{Initialization for nonvicously damped oscillators}
The integral constitutive relations of viscoelastic materials are represented by the following integro-differential equation of Volterra type:
\begin{equation}
\sigma \left( t \right)=\int_{-\infty }^{t}{G\left( t-\tau  \right)\dot{\varepsilon }\left( \tau  \right)d\tau }
\end{equation}
where $\sigma \left( t \right)$ is the stress, $\varepsilon \left( t \right)$ is the stain, $G\left( t \right)$ is the stress relaxation function. The lower terminal in the integral is $-\infty$ because the stress of viscoelastic materials is dependent on all the time histories of the strain [2].
When the integral constitutive equation (1) is applied to model the dynamics of structures incorporated with viscoelastic dampers, the equation of motion is
\begin{equation}
m\ddot{x}(t)+c\int_{-\infty }^{t}{G\left( t-\tau  \right)\dot{x}\left( \tau  \right)d}\tau +kx(t)=f\left( t \right)
\end{equation}
where $m$is the mass, $k$ is the stiffness, $c$ is the damping coefficient, $f\left( t \right)$ is the external force acting on the system.

The integral term in Eq. (2) makes the dynamic models different from the classical ones. It contains not only the information of vibrating displacement $x\left( t \right)$ and velocity $\dot{x}\left( t \right)$, but also the time-histories of velocity $\dot{x}\left( t \right)$. This implies that, unlike the viscously damped systems, the equation of motion, the instantaneous displacement and velocity are insufficient to predict the dynamic behaviors. Time-histories of motion should be added to initial conditions to fully determine the dynamics of nonvicously damped oscillators. As a result, the dynamic equation with past history is described as
\begin{equation}
\begin{cases}
m\ddot{x}(t)+c\int_{-\infty }^{t}{G\left( t-\tau  \right)\dot{x}\left( \tau  \right)d}\tau +kx(t)=f\left( t \right),\text{  }t>0\\
x(0)={{x}_{0}},\text{  }\dot{x}(0)={{v}_{0}}\\
\dot{x}(t)=v(t),\text{  }-\infty \ <t<0
\end{cases}
\end{equation}
where $t=0$ is the initial time and the lower terminal in the integral $t=-\infty $ is the starting time of vibration. In reality it is more reasonable to set the starting time as $t=-a$, which means that the system is at quiescent before $t=-a$ and begins vibrating at $t=-a$.

To be more familiar with the dynamic systems with memory, the integral term in Eq.(3) can be separated into two part:
\[\int_{-a}^{t}{G\left( t-\tau  \right)\dot{x}\left( \tau  \right)d}\tau =\int_{-a}^{0}{G\left( t-\tau  \right)\dot{x}\left( \tau  \right)d}\tau +\int_{0}^{t}{G\left( t-\tau  \right)\dot{x}\left( \tau  \right)d}\tau \]
The first part characterizes the hereditary effects of the histories of motion on the system dynamics, which is denoted as $\psi \left( t \right)$:
\begin{equation}
\psi \left( t \right)=\int_{-a}^{0}{G\left( t-\tau  \right)\dot{x}\left( \tau  \right)d}\tau =\int_{-a}^{0}{G\left( t-\tau  \right)v\left( \tau  \right)d}\tau
\end{equation}

It contains the time histories of motion, acts as an internal force and influence the behavior of dynamic systems after the initial time $t=0.$
Eq.(3) can be rewritten as
\begin{equation}\begin{cases}
m\ddot{x}(t)+c\int_{0}^{t}{G\left( t-\tau  \right)\dot{x}\left( \tau  \right)d}\tau +kx(t)=f\left( t \right)-c\psi \left( t \right),\text{  }t>0\\
x(0)={{x}_{0}},\text{  }\dot{x}(0)={{v}_{0}}
\end{cases}
\end{equation}
\newdefinition{rmka}{Remark}
\begin{rmka}
In studies of Adhikari and his colleagues [6], vibrating motions from past histories have not been taken into account. In this case, the equation of motion is
\[m\ddot{x}(t)+c\int_{0}^{t}{G\left( t-\tau  \right)\dot{x}\left( \tau  \right)d}\tau +kx(t)=f\left( t \right)\]
The initial values contain only the initial displacement $x(0)$ and the initial velocity $\dot{x}(0)$.
\end{rmka}
\begin{rmka}
In the integral constitutive relations of viscoelastic materials (1) and the corresponding dynamic equation of viscoelastically/nonviscously damped oscillators (2), $G\left( t \right)$ can be fitted by many types of decaying functions, such as the exponential functions, the power-law functions, the Mittage-Leffler functions, and so on. In the following sections, we concentrate on the case of exponentially decaying functions $G\left( t \right)=\mu {{e}^{-\mu t}},\mu >0$. The corresponding oscillators are called exponentially damped oscillators.
\end{rmka}
\section{Hereditary effects on the dynamic response}
Now we are ready to study the initialization response of exponentially damped oscillators with memories, which characterizes the hereditary effects of the vibrating motion from the past history. For this purpose, we will not take into account the external acting force $f\left( t \right)$ and set it to be zero. In this case, the equation of motion with initial condition is
\begin{equation}\begin{cases}
 m\ddot{x}(t)+c\int_{0}^{t}{G\left( t-\tau  \right)\dot{x}\left( \tau  \right)d}\tau +kx(t)=-c\psi \left( t \right),\text{  }t>0 \\
x(0)={{x}_{0}},\text{  }\dot{x}(0)={{v}_{0}}
\end{cases}
\end{equation}
where $G\left( t \right)=\mu {{e}^{-\mu t}},\mu >0$.
Taking Laplace transform of Eq.(6), one derives
\begin{equation}
m\left( {{s}^{2}}\bar{x}\left( s \right)-s{{x}_{0}}-{{v}_{0}} \right)+c\left( \frac{\mu s}{\mu +s}\bar{x}\left( s \right)-\frac{\mu }{\mu +s}{{x}_{0}} \right)+k\bar{x}\left( s \right)=-c\bar{\psi }\left( s \right)
\end{equation}
where $\bar{x}\left( s \right)$ is the Laplace transform of $x\left( t \right)$, $\bar{\psi }\left( s \right)$ is the Laplace transform of $\psi \left( t \right)$.
After rearranging Eq.(7), one has
\begin{equation}
\left( m{{s}^{2}}+\frac{c\mu s}{\mu +s}+k \right)\bar{x}\left( s \right)=-c\bar{\psi }\left( s \right)+ms{{x}_{0}}+\frac{c\mu }{\mu +s}{{x}_{0}}+m{{v}_{0}}
\end{equation}
We denote that $\bar{d}\left( s \right)=m{{s}^{2}}+\frac{c\mu s}{\mu +s}+k$ and $\bar{h}\left( s \right)=\frac{1}{\bar{d}\left( s \right)}$.\\
From Eq.(8), the solution of $\bar{x}\left( s \right)$ can be derived as
\begin{equation}
\bar{x}\left( s \right)=-c\bar{h}\left( s \right)\bar{\psi }\left( s \right)+m{{x}_{0}}s\bar{h}\left( s \right)+c\mu {{x}_{0}}\frac{\bar{h}\left( s \right)}{s+\mu }+m{{v}_{0}}\bar{h}\left( s \right)
\end{equation}
Taking the inverse Laplace transform of Eq.(9), one derives
\begin{equation}
\begin{split}
x\left( t \right)&=-c\int_{0}^{t}{h\left( t-\tau  \right)}\psi \left( \tau  \right)d\tau +m{{x}_{0}}\dot{h}\left( t \right)\\
& \quad +c\mu {{x}_{0}}\int_{0}^{t}{h\left( t-\tau  \right)}{{e}^{-\mu \tau }}d\tau +m{{v}_{0}}h\left( t \right)
\end{split}\end{equation}
where $h\left( t \right)$ is inverse Laplace transform of $\bar{h}\left( s \right)$.\\
Next we determine the expression of $h\left( t \right)$.
It has been shown in [6] that $\bar{d}\left( s \right)$ has zeros at $s={{s}_{j}},j=1,2,3$:
$${{s}_{1}}=-\alpha +\beta i,{{s}_{2}}=-\alpha -\beta i,{{s}_{1}}=-\gamma, $$
where $\alpha ,\beta ,\gamma >0$.\\
Furthermore, $\bar{h}\left( s \right)$can be expressed in the pole-residue form as
\[\bar{h}\left( s \right)\text{=}\sum\limits_{j=1}^{3}{\frac{{{R}_{j}}}{s-{{s}_{j}}}}\]
where ${{R}_{j}}$ are the residues and calculated as
\[{{R}_{i}}=\underset{s={{s}_{j}}}{\mathop{\text{Res}}}\,\bar{h}\left( s \right)=\underset{s\to {{s}_{j}}}{\mathop{\lim }}\,\left( s-{{s}_{j}} \right)\bar{h}\left( s \right)=\frac{1}{\underset{s\to {{s}_{j}}}{\mathop{\lim }}\,\frac{m{{s}^{2}}+\frac{c\mu s}{\mu +s}+k}{s-{{s}_{j}}}}=\frac{1}{\frac{\partial \bar{d}\left( s \right)}{\partial s}{{|}_{s={{s}_{j}}}}}\]
As a result, $h\left( t \right)$ can be obtained by taking Laplace transform of $\bar{h}\left( s \right)$:
\begin{equation}
h\left( t \right)={{L}^{-1}}\left\{ \bar{h}\left( s \right) \right\}\text{=}\sum\limits_{j=1}^{3}{{{R}_{j}}{{e}^{-{{s}_{j}}t}}}
\end{equation}
Substituting Eq.(11) into Eq.(10), one derives the initialization response of exponentially damped oscillators. It represents the hereditary effects of past histories of motions from the starting time at $t=-a$.
\section{Stability of initialization response }
It has been shown in Section 2 and Section 3 the effects of past histories of motions on the setting of initial conditions and dynamic response. In this section, we proceed to show the hereditary effects on stability of initialization response. We prove that this influence will gradually recede with increasing of time.
From Eq.(10), it is clear that the initialization response involves four parts. The last three parts decreased as time increases. Now we prove in mathematics that the same holds for the first part.
\begin{equation}
\begin{split}
&\quad-c\int_{0}^{t}{h\left( t-\tau  \right)\psi \left( \tau  \right)d\tau}\\
 &=-c\int_{0}^{t}{{{R}_{1}}{{e}^{{{s}_{1}}\left( t-\tau  \right)}}\psi \left( \tau  \right)d\tau }-c\int_{0}^{t}{{{R}_{1}}{{e}^{{{s}_{2}}\left( t-\tau  \right)}}\psi \left( \tau  \right)d\tau }\\
 &\quad-c\int_{0}^{t}{{{R}_{3}}{{e}^{{{s}_{3}}\left( t-\tau  \right)}}\psi \left( \tau  \right)d\tau } \\
 &=-c\int_{0}^{t}{{{R}_{1}}{{e}^{-\left( \alpha -\beta i \right)\left( t-\tau  \right)}}\psi \left( \tau  \right)d\tau }-c\int_{0}^{t}{{{R}_{1}}{{e}^{-\left( \alpha +\beta i \right)\left( t-\tau  \right)}}\psi \left( \tau  \right)d\tau }\\
 &\quad-c\int_{0}^{t}{{{R}_{3}}{{e}^{-\gamma \left( t-\tau  \right)}}\psi \left( \tau  \right)d\tau } \\
 &=-2c{{R}_{1}}\int_{0}^{t}{{{e}^{-\alpha \left( t-\tau  \right)}}\cos \beta \left( t-\tau  \right)\psi \left( \tau  \right)d\tau }-c{{R}_{3}}\int_{0}^{t}{{{e}^{-\gamma \left( t-\tau  \right)}}\psi \left( \tau  \right)d\tau } \\
 &={{I}_{1}}+{{I}_{2}}
\end{split}
\end{equation}
where
\begin{equation}
{{I}_{1}}=-2c{{R}_{1}}\int_{0}^{t}{{{e}^{-\alpha \left( t-\tau  \right)}}\cos \beta \left( t-\tau  \right)\psi \left( \tau  \right)d\tau }
\end{equation}
\begin{equation}
{{I}_{2}}=-c{{R}_{3}}\int_{0}^{t}{{{e}^{-\gamma \left( t-\tau  \right)}}\psi \left( \tau  \right)d\tau }
\end{equation}
Substituting Eq.(4) into Eq.(13),one has
\begin{equation}
\begin{split}
\left| {{I}_{1}} \right|&=2c{{R}_{1}}\left| \int_{0}^{t}{{{e}^{-\alpha \left( t-\tau  \right)}}\cos \beta \left( t-\tau  \right)\psi \left( \tau  \right)d\tau } \right| \\
 &=2c{{R}_{1}}\left| \int_{0}^{t}{{{e}^{-\alpha \left( t-\tau  \right)}}\cos \beta \left( t-\tau  \right)d\tau \int_{-a}^{0}{G\left( \tau -{{\tau }_{1}} \right)}v\left( {{\tau }_{1}} \right)d{{\tau }_{1}}} \right| \\
\end{split}\end{equation}

It is reasonable to suppose that the response velocity before initial time is bounded, i.e.,
$\left| v\left( t \right) \right|\le M,t\in \left[ -a,0 \right]$.

Noting that$\left| \cos \left( \beta t \right) \right|\le 1$, Eq.(15) yields
\begin{equation}
\left| {{I}_{1}} \right|\le 2c{{R}_{1}}M\left| \int_{0}^{t}{{{e}^{-\alpha \left( t-\tau  \right)}}d\tau \int_{-a}^{0}{G\left( \tau -{{\tau }_{1}} \right)}d{{\tau }_{1}}} \right|
\end{equation}

Substituting $G\left( t \right)=\mu {{e}^{-\mu t}}$ into Eq.(16) yields
\begin{equation}\begin{split}
\left| {{I}_{1}} \right|&\le 2c{{R}_{1}}M\left( 1-{{e}^{-\mu a}} \right)\left| \int_{0}^{t}{{{e}^{-\alpha \left( t-\tau  \right)-\mu \tau }}d\tau } \right| \\
 & =2c{{R}_{1}}M\left( 1-{{e}^{-\mu a}} \right){{e}^{-\alpha t}}\left| \int_{0}^{t}{{{e}^{\left( \alpha -\mu  \right)\tau }}d\tau } \right| \\
 & =\frac{2c{{R}_{1}}M\left( 1-{{e}^{-\mu a}} \right)}{\left| \alpha -\mu  \right|}\left| {{e}^{-\mu t}}-{{e}^{-\alpha t}} \right| \\
\end{split}\end{equation}

Due to the fact that $\mu ,\alpha >0$, it is clear to see that
\begin{equation}
\underset{t\to \infty }{\mathop{\lim }}\,\left| {{I}_{1}} \right|=0.
\end{equation}
Substituting Eq.(4) into Eq.(14), one has
\begin{equation}\begin{split}
\left| {{I}_{2}} \right|&=c{{R}_{3}}\left| \int_{0}^{t}{{{e}^{-\gamma \left( t-\tau  \right)}}\psi \left( \tau  \right)d\tau } \right| \\
 & \text{=}c{{R}_{3}}\left| \int_{0}^{t}{{{e}^{-\gamma \left( t-\tau  \right)}}d\tau \int_{-a}^{0}{G\left( \tau -{{\tau }_{1}} \right)}v\left( {{\tau }_{1}} \right)d{{\tau }_{1}}} \right| \\
 & \le Mc{{R}_{3}}\left( 1-{{e}^{-\mu a}} \right)\left| \int_{0}^{t}{{{e}^{-\gamma \left( t-\tau  \right)-\mu \tau }}d\tau } \right| \\
 & =Mc{{R}_{3}}\left( 1-{{e}^{-\mu a}} \right){{e}^{-\gamma t}}\left| \int_{0}^{t}{{{e}^{\left( \gamma -\mu  \right)\tau }}d\tau } \right| \\
 & =\frac{Mc{{R}_{3}}\left( 1-{{e}^{-\mu a}} \right)}{\left| \gamma -\mu  \right|}\left| {{e}^{-\mu t}}-{{e}^{-\gamma t}} \right| \\
\end{split}\end{equation}
Because $\mu ,\gamma >0$, it is also clear to see that
\begin{equation}
\underset{t\to \infty }{\mathop{\lim }}\,\left| {{I}_{2}} \right|=0
\end{equation}
From (12) (18) and (20), we derive that the first part of the initialization response (10) decreases to zero. It is very easy to verify that the remaining three parts of Eq.(10) also decreases to zero. So, we can conclude that the initialization response will gradually recede with increasing of time. In other words, memories from past histories of motion have no influence on the stability of exponentially damped oscillators.

\section{Conclusions}
Hereditary effects of exponentially damped oscillators with past histories have been investigated in this paper. Initial conditions have been proposed for nonvicously damped oscillators with past histories. It has been shown that initial conditions should contain not only the vibration displacement and velocity at initial time, but also the time-history of response velocity from the starting time of vibration. Then, initialization response of exponentially damped oscillators has been obtained to characterize the hereditary effects of past histories on the dynamic response. At last, initialization response has been proved to gradually recede with increasing of time, which implies that memories from past histories have no influence on the stability of exponentially damped oscillators.\\
\textbf{Acknowledgements}
All the authors acknowledge the valuable suggestions from the peer reviewers. This work was supported by the National Natural Science Foundation of China (Grant No. 11802338).


\begin{thebibliography}{99}
\bibitem{1}Ibrahim, R.A., Recent advances in nonlinear passive vibration isolators, Journal of Sound Vibration 314 (2008) 371-452.
\bibitem{2}Shaw, M.T. and W.J. Macknight, Introduction to polymer viscoelasticity, Wiley, 2005.
\bibitem{3}Woodhouse, Linear damping models for structural vibration, Journal of Sound Vibration 215(1998) 547-569.
\bibitem{4}Muravyov, Forced vibration responses of a viscoelastic structure, Journal of Sound Vibration 218(1998) 892-907.
\bibitem{5}Mainardi, F., Fractional Calculus and Waves in Linear Viscoelasticity, Imperial College Press, 2010.
\bibitem{6}Adhikari, S., Structural Dynamic Analysis with Generalized Damping Models, Wiley, 2014.
\bibitem{7}Adhikari, S., Dynamic Response Characteristics of a Nonviscously Damped Oscillator, Journal of Applied Mechanics 75(2008) 148-155.
\bibitem{8}Adhikari, S. and J. Woodhouse, Quantification of non-viscous damping in discrete linear systems, Journal of Sound Vibration, 260(2003) 499-518.
\bibitem{9}Garc\'{\i}a-Barruetabe\~{n}a, J., et al., Dynamics of an exponentially damped solid rod: Analytic solution and finite element formulations, International Journal of Solids \& Structures, 49(2012) 590-598.
\bibitem{10}Li L, Hu Y J, Wang X L, et al. Computation of Eigensolution Derivatives for Nonviscously Damped Systems Using the Algebraic Method, AIAA Journal, 50(2012) 2282-2284.
\bibitem{11}Reggio, A., et al., A state-space methodology to identify modal and physical parameters of non-viscously damped systems, 41(2013): p. 380-395.
\bibitem{12}L\'{a}zaro, M., Closed-form eigensolutions of nonviscously, nonproportionally damped systems based on continuous damping sensitivity. Journal of Sound Vibration, 413(2018) 368-382.
\bibitem{13}Rossikhin Y A, Shitikova M V, Application of Fractional Calculus for Dynamic Problems of Solid Mechanics: Novel Trends and Recent Results,Applied Mechanics Reviews, 2009, 63(1):010801.
\bibitem{14}Paola, M.D., A. Pirrotta, and A.J.M.o.M. Valenza, Visco-elastic behavior through fractional calculus: An easier method for best fitting experimental results. 2011. 43(12): p. 799-806.
\bibitem{15}Rossikhin, Y.A. and M.V. Shitikova, Application of fractional calculus for dynamic problems of solid mechanics: novel trends and recent results. Applied Mechanics Reviews, 2010. 63(1): p. 010801.
\bibitem{16}Padovan, J., S. Chung, and Y.H. Guo, Asymptotic steady state behavior of fractionally damped systems. Journal of the Franklin Institute, 1987. 324(3): p. 491-511.
\bibitem{17}Padovan, J. and Y. Guo, General response of viscoelastic systems modelled by fractional operators. Journal of the Franklin Institute, 1988. 325(2): p. 247-275.
\bibitem{18}	Beyer, H. and S. Kempfle, Definition of physically consistent damping laws with fractional derivatives. ZAMM Journal of Applied Mathematics and Mechanics/Zeitschrift f\"{u}r Angewandte Mathematik und Mechanik, 1995. 75(8): p. 623-635.
\bibitem{19}Kempfle, S., I. Sch\"{a}fer, and H. Beyer, Fractional calculus via functional calculus: theory and applications. Nonlinear Dynamics, 2002. 29(1-4): p. 99-127.
\bibitem{20}Sch\"{a}fer, I. and S. Kempfle, Impulse responses of fractional damped systems. Nonlinear Dynamics, 2004. 38(1-4): p. 61-68.
\bibitem{21}Rossikhin, Y.A. and M.V. Shitikova, Analysis of rheological equations involving more than one fractional parameters by the use of the simplest mechanical systems based on these equations. Mechanics of Time-Dependent Materials, 2001. 5(2): p. 131-175.
\bibitem{22}Yuan J, Zhang Y, Liu J, et al. Mechanical energy and equivalent differential equations of motion for single-degree-of-freedom fractional oscillators. 2017. 397: p. 192-203.
\bibitem{23}Zhang Y, Yuan J, Liu J, et al. Lyapunov Functions and Sliding Mode Control for Two Degrees-of-Freedom and Multidegrees-of-Freedom Fractional Oscillators. 2017. 139(1).
\bibitem{24}Yuan J,  Sliding Mode Control of Vibration in Single-Degree-of-Freedom Fractional Oscillators. 2017. 139(11): p. 114503-114503-6.
\bibitem{25}Podlubny, I. Fractional Differential Equations. in Mathematics in Science and Engineering. 1999.
\bibitem{26}Fukunaga M. On initial value problems of fractional differential equations. International Journal of Applied Mathematics, 2002, 9(2): 219-236
\bibitem{27}Fukunaga M. On uniqueness of the solutions of initial value problems of ordinary fractional differential equations. International Journal of Applied Mathematics, 2002, 10(2): 177-190
\bibitem{28}Fukunaga M. A difference method for initial value problems for ordinary fractional differential equations, II. International Journal of Applied Mathematics, 2003, 11(3): 215-244
\bibitem{29}Hartley, T.T. and C.F. Lorenzo, Control of initialized fractional-order systems. NASA Technical Report, 2002.
\bibitem{30}Lorenzo, C.F. and T.T. Hartley, Initialization of Fractional-Order Operators and Fractional Differential Equations. Journal of Computational and Nonlinear Dynamics, 2008. 3(2): p. 021101.
\bibitem{31}Wu C, Yuan J, Shi B. Stability of initialization response of fractional oscillators, Journal of Vibroengineering, 2016.
\end{thebibliography}
\end{document}